\documentclass[11pt, a4paper]{article}
\usepackage{amsmath,amssymb,latexsym,graphics,epsfig}
\usepackage{hyperref}
\usepackage{color}
\usepackage{amsthm}

\newtheorem{theo}{Theorem}[section]
\newtheorem{propo}[theo]{Proposition}

\newtheorem{coro}[theo]{Corollary}

\input amssym.def
\newsymbol\rtimes 226F
\newfont{\nset}{msbm10}
\newcommand{\ns}[1]{\mbox{\nset #1}}

\def\Z{\ns Z}

\def\A{\mbox{\boldmath $A$}}
\def\B{\mbox{\boldmath $B$}}

\def\H{\mbox{\boldmath $H$}}
\def\I{\mbox{\boldmath $I$}}

\def\M{\mbox{\boldmath $M$}}

\def\U{\mbox{\boldmath $U$}}
\def\V{\mbox{\boldmath $V$}}

\def\S{\mbox{\boldmath $S$}}

\def\X{\mbox{\boldmath $X$}}
\def\Z{\ns{Z}}

\def\>{\mathop{\rightarrow}\nolimits}
\def\a{\mbox{\boldmath $a$}}
\def\b{\mbox{\boldmath $b$}}
\def\e{\mbox{\boldmath $e$}}

\def\g{\mbox{\boldmath $g$}}
\def\i{\mbox{\boldmath $i$}}

\def\m{\mbox{\boldmath $m$}}

\def\u{{\mbox {\boldmath $u$}}}
\def\v{\mbox{\boldmath $v$}}
\def\w{{\mbox {\boldmath $w$}}}
\def\x{\mbox{\boldmath $x$}}

\def\z{\mbox{\boldmath $z$}}

\def\vec0{\mbox{\bf 0}}

\def\Aut{\mathop{\rm Aut}\nolimits}

\begin{document}

\title{On Congruence in $Z^n$ and the Dimension of a
Multidimensional Circulant}
\author{M.A. Fiol
\\ \\
{\small Universitat Polit\`ecnica de Catalunya, BarcelonaTech}\\
{\small Departament de Matem\`atica Aplicada IV} \\
{\small Barcelona, Catalonia} \\
{\small (e-mail: {\tt fiol@ma4.upc.edu})}}

\date{}

\maketitle

\begin{abstract}

From a generalization to $\Z^n$ of the concept of congruence we define
a family of regular digraphs or graphs called multidimensional
circulants, which turn out to be Cayley (di)graphs of Abelian groups.
This paper is mainly devoted to show the relationship between the
Smith normal form for integral matrices and the dimensions of such
(di)graphs, that is the minimum ranks of the groups they can arise from.
In particular, those $2$-step multidimensional circulants which are
circulants, that is Cayley (di)graphs of cyclic groups, are fully
characterized. In addition, a reasoning due to Lawrence is used to
prove that the cartesian product of $n$ circulants with equal number of
vertices $p>2$, $p$ a prime, has dimension~$n$.

\end{abstract}


\section{Introduction}

Throughout this paper we make use of standard concepts and
terminology of graph theory and group theory, see for instance \cite{r6}
and \cite{r15} respectively. We will recall here the most
relevant definitions. Let $\Gamma$ be a group with identity element $e$,
and let $A\subseteq \Gamma\setminus \{ e\}$ such that $A^{-1} = A$.
The Cayley graph of $\Gamma$ with respect to $A$, denoted by
$G(\Gamma ;A)$, is the graph whose vertices are labelled with the
elements of $\Gamma$, and an edge $(u,v)$ if and only if $u^{-1}v \in
A$. The Cayley digraph $G(\Gamma ;A)$ is defined similarly, but now we
do not require $A^{-1} = A$. Since left translations in $\Gamma$ are
automorphisms of $G(\Gamma ;A)$, a Cayley graph is always
vertex-transitive. Moreover, the group of such automorphisms,
$\Aut G(\Gamma ;A)$, contains a regular subgroup (that is a
transitive subgroup whose order coincides with the order of the graph)
isomorphic to $\Gamma$. In fact, Sabidussi \cite{r16} first showed that
this is also a sufficient condition for a graph to be a Cayley graph
(of the group $\Gamma$). Clearly, the same statements also hold for
Cayley digraphs. In particular, this result or its digraph analog
will be simply referred to as {\it Sabidussi's result}.

Because of both theoretical and practical reasons, the class of
Cayley (di)\-graphs obtained from Abelian groups have deserved special
attention in the literature. This is the case, for instance, for
circulant (di)graphs the definition of which follows. Let $m$ be a
positive integer, and $A =\{ a_{1} ,a_{2} , \ldots ,a_{d} \}
\subseteq \Z/m\Z\setminus \{ 0\}$. The {\it $($d-step$)$ circulant digraph} $G(m;A)$
has as set of vertices the integers modulo $m$, and vertex $u$ is
adjacent to the vertices $u+A = \{ u+a_i\pmod{m} : a_i \in A \}$. The names
{\it multiple loop digraph} and {\it multiple fixed step digraph}
\cite{r3,r12} are also used. The {\it $($d-step$)$ circulant graph} $G(m;A)$ is
defined analogously with $A = -A = \{ \pm a_{1}, \pm a_{2}, \ldots ,
\pm a_{d} \}$. These graphs have also received other names
such as {\it starred polygons} \cite{r17} and {\it multiple loop graphs}
\cite{r4}. In both the directed and undirected case, the elements of
$A$ are called {\it jumps} or {\it steps}. Henceforth, we will use the
word {\it circulant} to denote either a circulant digraph or a circulant
graph. As stated above, notice that the circulant $G(m;A)$ is just the
Cayley (di)graph of the cyclic group $Z/mZ$ with respect to $A$.
Therefore, from Sabidussi's result,
{\it a $($di$)$graph is a circulant iff its automorphism group contains a
regular cyclic subgroup}. As Leighton showed in \cite{r13}, this
characterization can be used to easily prove Turner's result \cite{r17},
which states that every transitive graph on a prime number of
vertices is a circulant. (It suffices to use Cauchy's group
theorem: `If a prime $p$ divides the order of a finite group then it
contains an element of order $p$.') In fact the same reasoning shows
that Turner's result also holds for digraphs.

Recall that, given two graphs, $G_{1} = (V_{1},E_{1})$ and $G_{2} =
(V_{2},E_{2})$, their {\it cartesian product} $G_{1} \times G_{2}$ is the
graph with set of vertices $V_{1} \times V_{2}$ and an edge between
$(u_{1},u_{2})$ and $(v_{1},v_{2})$ iff either $(u_{1},v_{1})
\in E_{1}$ and $u_{2} =v_{2}$ or $u_{1} = v_{1}$ and $(u_{2},v_{2})
\in E_{2}$. The cartesian product of two digraphs is defined
analogously. Another well known family of Cayley graphs of Abelian
groups are the (Boolean) $n$-cubes, which are sometimes defined as the
cartesian product of $n$ copies of the complete graph $K_{2}$. From our
point of view, the {\it n-cube} (or {\it binary n-dimensional
hypercube}) is the Cayley graph $G(\Gamma ;A)$, where
$\Gamma = \Z/2\Z \times \cdots \times \Z/2\Z$ ($n$ factors) and $A$ is
the set of unitary vectors $\e_{i}$, $1\leq i \leq n$.

This paper studies the structure of the so-called ($d$-step)
multidimensional circulants (that is Cayley graphs or digraphs of
Abelian groups), which are defined from integral matrices in Section 3.
More precisely, given such a (di)graph, we are interested in finding
its `dimension', that is the minimum rank of the group(s) it can
arise from. (The {\it rank} of a finitely generated Abelian group is
the minimum number of elements which generate it.) For instance, from
its definition it is clear that the $n$-cube has dimension not greater
than $n$. As far as we know, this study was initialized by Leighton
\cite{r13}, where it was shown that the dimension of the $n$-cube is in
fact $\lfloor (n+1)/2 \rfloor$. Other results, concerning $2$-step
circulant digraphs (dimension $1$) can be found in \cite{eaf,r12}. In this
paper we continue this study by deriving some new results, which are
contained in Section 3. For instance, using some facts about integral
matrices and a theorem in \cite{r7}, we give a full characterization of
those $2$-step multidimensional circulants which are $1$-dimensional.
Moreover, it is proved that if $G_{1}, G_{2}, \ldots ,G_{n}$ are
circulant (di)graphs with $p>2$ vertices, $p$ a prime, then the
cartesian product $G_{1}\times G_{2}\times \cdots \times G_{n}$ has
dimension $n$.

As stated above, our approach uses some results of integral matrix
theory which are summarized in the next section. In particular, we deal
with the concept of congruence in $\Z^{n}$, also discussed there. The
reason is that, in the same way that congruence in $\Z$ (periodicity in one
dimension) leads to the consideration of cyclic groups, congruence in
$\Z^{n}$ (related to periodicity in $n$ dimensions) induces quotient
structures which are Abelian groups.


\section{Congruences in $Z^n$ and the
Induced Abelian Groups}

This section mainly deals with the concept of congruence in $\Z^{n}$
and its consequences to our study. In this context, the theory of
integral matrices (that is matrices whose entries are integers) proves
to be very useful and, hence, we begin by recalling some of its main
results. The reader interested in the proofs is referred to \cite{r14}.

\medskip
Let $\Z^{n\ast n}$ be the ring of $n \times n$ matrices over $\Z$. Given
$\A,\B\in \Z^{n\ast n}$, $\A$ is said to be {\it right equivalent} to $\B$
if there exists a unimodular (with determinant $\pm 1$) matrix $\V\in
\Z^{n\ast n}$ such that $\A=\B\V$; and $\A$ is {\it equivalent} to $\B$ if
$\A=\U\B\V$ for some unimodular matrices $\U,\V\in \Z^{n\ast n}$. Clearly,
both of them are equivalence relations.

Henceforth, $\M=(m_{ij})$ will denote a nonsingular matrix of
$\Z^{n\ast n}$ with columns $\m_{j}=(m_{1j},m_{2j},\ldots
,m_{nj})^\top$, $j=1,2, \ldots ,n$, and $m=|\det \M|$. By the Hermite
normal form theorem, $\M$ is right equivalent to an upper triangular
matrix $\H(\M)=\H=(h_{ij})$ with positive diagonal elements $h_{ii}$ and
with each element above the main diagonal $h_{ij}$,  $j>i$, $i=1,2,\ldots
,n-1$, lying in a given complete set of residues modulo $h_{ii}$ (for
instance, $0\leq h_{ij} \leq h_{ii}-1)$. This normal form is unique.

Let $k\in \Z$, $1\leq k\leq n$. The $k$th {\it determinantal divisor} of
$\M$, denoted by $d_{k}(\M)=d_{k}$, is defined as the greatest common
divisor of the $(^{n}_{k})^{2}$  $k\times k$ determinantal minors of $\M$.
Since $\M$ is nonsingular, not all of them are zero. Notice that $d_{k}|
d_{k+1}$ for all $k=1,2,\ldots ,n-1$ and $d_{n}=m$. For convenience,
put $d_{0}=1$. The {\it invariant factors} of $\M$ are the quantities
$$
s_{k}(\M) = s_{k} = \frac{d_{k}}{d_{k-1}},\quad k=1,2,\ldots ,n.
$$
It can be shown that $s_{i}|s_{i+1}$, $i=1,2,\ldots ,n-1$.

By the Smith normal form theorem, $\M$ is equivalent to the diagonal
matrix $\S(\M)=\S=\;$diag$(s_{1},s_{2},\ldots ,s_{n})$. This canonical form
is unique.

For instance, the matrix $\M=\;$diag$(2,2,3)$, with determinantal
divisors $d_{1}=1$, $d_{2} =2$, $d_{3}=12$, and invariant factors $s_{1} =1$,
$s_{2}=2$, $s_{3}=6$, is equivalent to $\S=\;$diag$(1,2,6)$ since there exist
the unimodular matrices
\begin{equation}
\U=\left( \begin{array}{rrr}
-1 & 0 & 1 \\ 0 & 1 & 0 \\ -3 & 0 & 2
         \end{array} \right),\quad
\V=\left( \begin{array}{rrr}
1 & 0 & 3 \\ 0 & 1 & 0 \\ 1 & 0 & 2
         \end{array} \right)                                                       \label{eq1}
\end{equation}
such that $\S=\U\M\V$.

As usual, the greatest common divisor of the integers $a_{1},a_{2},
\ldots ,a_{n}$ will be denoted by $\gcd(a_{1},a_{2}, \ldots ,a_{n} )$.
When they are the coordinates of a vector $\a$, we will simply
write $\gcd(\a)$.

Most of the remaining material in this section is drawn from \cite{r9}
where additional details can be found.

Let $\Z^{n}$ denote the additive group of column $n$-vectors with
integral coordinates. The set $\M\Z^{n}$ , whose elements are linear
combinations (with integral coefficients) of the (column) vectors
$\m_{j}$, is said to be the lattice generated by $\M$. Clearly,
$\M\Z^{n}$ with the usual vector addition is a normal subgroup of $\Z^{n}$.

The concept of congruence in $\Z$ has the following natural
generalization to $\Z^{n}$. Let $\a,\b\in \Z^{n}$ .We say that
{\it $\a$  is congruent with $\b$ modulo $\M$}, and write
$\a \equiv \b\pmod{\M}$, if
\begin{equation}
\label{eq2}
\a-\b \in \M\Z^{n}.
\end{equation}
The quotient group $\Z^{n}/\M\Z^{n}$ can intuitively be called the
{\it group of integral vectors modulo \M}. Henceforth, we follow the
usual convention of identifying each equivalence class by any of its
representatives.

Note that whenever $\M=\;$diag$(m_{1},m_{2}, \ldots ,m_{n})$ the vectors
\a=$(a_{1},a_{2}, \ldots ,a_{n})^\top$ and
\b=$(b_{1},b_{2}, \ldots ,b_{n})^\top$
are congruent modulo $\M$ iff the system of congruences in $\Z$
$$
 a_{i} \equiv b_{i}\pmod{m_{i}},\quad i=1,2, \ldots ,n
 $$
holds. In this case $\Z^{n}/\M\Z^{n}$ is the direct product of the cyclic
groups $\Z/m_{i}\Z$, $i=1,2, \ldots ,n$.

If $\A$ and $\B$ are $n\times r$ matrices over $\Z$ with columns
$\a_{j}$ and $\b_{j}$, $j=1,2,\ldots ,r$, respectively, we will
write $\A \equiv \B\pmod{\M}$ to denote that
$\a_{j} \equiv \b_{j}\pmod{\M}$ for all $j=1,2,\ldots ,r$.

Let $\H=\M\V$ be the Hermite normal form of $\M$. Then (\ref{eq2}) holds
iff $\a-\b \in \H\V^{-1}\Z^{n} = \H\Z^{n}$ since $\V$, and hence $\V^{-1}$,
are unimodular. Therefore we conclude that
\begin{equation}
\a \equiv \b\pmod{\M}\quad \Leftrightarrow \quad
\a \equiv \b\pmod{\H}                                             \label{eq3}
\end{equation}
or, what is the same,
\begin{equation}
\Z^{n}/\M\Z^{n} \cong \Z^{n}/\H\Z^{n}.                                              \label{eq4}
\end{equation}

Let us now consider the Smith normal form of $\M$, $\S=\;$diag$(s_{1},s_{2},
\ldots ,s_{n})=\U\M\V$. Then (\ref{eq2}) holds iff $\U\a \equiv
\U\b\pmod{\S}$ or, equivalently,
\begin{equation}
\u_{i}\a \equiv \u_{i}\b\pmod{s_i},\quad
i=1,2, \ldots ,n                                                                      \label{eq5}
\end{equation}
where $\u_{i}$ stands for the $i$th row of $\U$. Moreover, if $r$
is the smallest integer such that $s_{n-r}=1$ (thus, $s_{1}=s_{2}=
\cdots =s_{n-r-1} = 1$), (if there is no such a $r$, let $r=n$), the
first $n-r$ equations in (\ref{eq5}) are irrelevant, and we only need to
consider the other ones. This allows us to write
\begin{equation}
\a \equiv \b\pmod{\M} \quad \Leftrightarrow \quad
\U'\a \equiv \U'\b\pmod{\S'}                                         \label{eq6}
\end{equation}
where $\U'$ stands for the $r\times n$ matrix obtained from $\U$ by
leaving out the first $n-r$ rows, and $\S'=\;$diag$(s_{n-r+1},s_{n-r+2},
\ldots ,s_{n})$. So, the (linear) mapping $\phi$ from the vectors modulo
$\M$ to the vectors modulo $\S'$ given by $\phi (\a)=\U'\a$ is a
group isomorphism, and we can write
\begin{equation}
\Z^{n}/\M\Z^{n} \cong \Z^{r}/\S'\Z^{r} = \Z/s_{n-r+1}\Z \times \cdots \times \Z/s_{n}\Z.
                                                                                                 \label{eq7}
\end{equation}
Analogously, it may be shown that the $n \times r$ matrix of the inverse
mapping $\phi^{-1}$ is obtained from $\U^{-1}$ by leaving out its first
$n-r$ columns.

The next proposition contains some easy consequences of the above
results. For instance, $(b)$ follows from the fact that
$s_{1}s_{2} \cdots s_{n} = d_{n} = m$ and $s_{i}|s_{i+1}$,
$i=1,2, \ldots ,n-1$.

\begin{propo}
\label{pro2.1}
\begin{itemize}
\item[$(a)$]
The number of equivalence classes modulo $\M$ is
$ |\Z^{n}/\M\Z^{n}| = m = |\det \M| $.
\item[$(b)$]
 If $p_{1}^{r_{1}}p_{2}^{r_{2}} \cdots p_{t}^{r_{t}}$ is the prime
factorization of $m$, then $\Z^{n}/\M\Z^{n} \cong \Z^{r}/\S'\Z^{r}$ for some
$r\times r$ matrix $\S'$ with $r \leq \max \{r_{i} :1 \leq i \leq t \}$.
\item[$(c)$]
 The $($Abelian$)$ group of integral vectors modulo $\M$ is cyclic
iff $d_{n-1}=$1.
\item[$(d)$]
 Let $r$ be the smallest integer such that $s_{n-r}=1$. Then $r$ is
the rank of $\Z^{n}/\M\Z^{n}$ and the last $r$ columns of $\U^{-1}$ form a basis
of $\Z^{n}/\M\Z^{n}$.\ $\Box$
\end{itemize}
\end{propo}

Given any element \a of $\Z^{n}/\M\Z^{n}$, simple reasoning shows
that its order is given by the formula
\begin{equation}
 \label{eq8}
{\rm o}(\a) = \frac{m}{\gcd(m,\gcd(m\M^{-1}\a))},
\end{equation}
(see \cite{r9}). For instance, if $\M=(m_{ij})$ is a $2\times 2$ matrix
and $\a=(a_{1},a_{2})^\top$ we have
\begin{equation}
{\rm o}(\a) = \frac{m}{\gcd(m,a_{1}m_{22}-a_{2}m_{12},
a_{2}m_{11}-a_{1}m_{21})}.                                                       \label{eq9}
\end{equation}
According to (\ref{eq7}), for any given $\M\in \Z^{n \ast n}$ there
exists an Abelian group $\Gamma$ such that $\Gamma=\Z^{n}/\M\Z^{n}$.
Conversely, let $\Gamma$ be a finite Abelian group generated by the
elements $g_{1},g_{2}, \ldots g_{n}$. Then $\Gamma$ is isomorphic to
$\Z^{n}/K$, where $K$ is the kernel of the surjective homomorphism
$\Psi : \Z^{n} \longrightarrow \Gamma$ defined by $\Psi(\x)=
\g\x$, where $\g$ denotes the row vector
$(g_{1},g_{2}, \ldots ,g_{n})$ and $\x\in \Z^{n}$. (Note that
$\Psi(\e_{i})=g_{i}$, $1\leq i\leq n$, where $\e_{i}$ stands for
the $i$th unitary vector.) More precisely, $K$ is the lattice of $\Z^{n}$
generated by the upper triangular $n\times n$ matrix $\H=(m_{ij})$
defined as follows:
\begin{itemize}
\item[]
 $m_{11}={\rm o}(g_{1})=|\langle g_{1}\rangle|$;
\item[]
 $m_{jj}=  \min\{  \mu  \in  \Z^{+}:\mu  g_{j}  \in  \langle g_{1},g_{2}, \ldots
,g_{j-1}\rangle \}$,  $j=2,3, \ldots ,n $; and
\item[]
$m_{ij}$, $j=2,3, \ldots ,n$, $i<j$, are any integers such that
$m_{1j}g_{1}+m_{2j}g_{2}+ \cdots +m_{jj}g_{j} = 0$, (they can be chosen
in a given complete set of residues modulo $h_{ii}$, e.g.,
$0 \leq m_{ij} \leq h_{ii} -1$),
\end{itemize}
where, as usual, $\langle g_{1},g_{2},
\ldots ,g_{j}\rangle$ denotes the group generated by $g_{1},g_{2}, \ldots
 ,g_{j}$, and $0$ is the identity element of $\Gamma$. Clearly, $\H$ is
the Hermite normal form of any matrix $\M$ which  generates  the  lattice
$K$.


\section{Multidimensional circulants and their dimension}

Congruence in $\Z^{n}$ leads to the following generalization of
circulants. Let $\M$ be an $n \times n$ integral matrix as in Section 2.
Let $A=\{ \a_{j}=(a_{1j},a_{2j}, \ldots ,a_{nj})^\top :
1 \leq j \leq d \} \subseteq \Z^{n}/\M\Z^{n}$. The {\it multidimensional
$($d-step$)$ circulant digraph G(\M;A)} has as vertex-set the integral
vectors modulo $\M$, and every vertex $\u$ is adjacent to the vertices
$\u+A\pmod{\M}$. As in the case of circulants, the {\it
multidimensional $($d-step$)$ circulant graph $G(\M;A)$} is defined similarly
just requiring $A=-A$.

In \cite{r13}, Leighton considered multidimensional circulant graphs with
diagonal matrix $\M$, and characterized them by showing that {\it the
automorphism group of these graphs must contain a regular Abelian
subgroup}. Clearly, a multidimensional circulant (digraph or graph) is a
Cayley (di)graph of the Abelian group $\Z^{n}/\M\Z^{n}$. As a consequence,
Sabidussi's result implies that Leighton's statement holds in fact for
a multidimensional circulant obtained from any matrix $\M$.

As another consequence of the above, if $\alpha$ is the index of the
subgroup \linebreak $\Gamma=\langle \a_{1},\a_{2}, \ldots ,\a_{d}\rangle$ in
$\Z^{n}/\M\Z^{n}$ , the multidimensional  circulant  $G(\M;A)$  consists  of
$\alpha$ copies of the Cayley (di)graph of $\Gamma$ generated by $A$.
Besides, from the comments in the last paragraph of Section 2,
$\Gamma \cong \Z^{d}/\H\Z^{d}$, where $\H$ is an upper triangular
$d \times d$ matrix, and each such copy is isomorphic to
$G(\H;\e_{1},\e_{2}, \ldots \e_{d})$.

In particular, $G(\M;\a_{1},\a_{2}, \ldots \a_{d})$
(respectively $G(\M;\pm\a_{1},\pm\a_{2}, \ldots \pm\a_{d})$)
is strongly connected (respectively connected), that is $\alpha=1$, iff
$\{ \a_{1},\a_{2}, \ldots \a_{d} \}$ generates
$\Z^{n}/\M\Z^{n}$, that is, there exist $n$ integral $d$-vectors
$\x^{j} = (x_{1j},x_{2j}, \ldots ,x_{dj})^\top$,
$j=1,2, \ldots ,n$, such that
$$
 x_{1j}\a_{1}+x_{2j}\a_{2}+ \cdots +x_{dj}\a_{d} \equiv
\e_{j}\pmod{\M},\quad j=1,2, \ldots ,n
$$
or, in matrix form,
$$
 \A\X \equiv \I\pmod{\M}
$$
where $\A$ now denotes the $n\times d$ matrix $(a_{ij})$, $\X$ is the
$d\times n$ matrix $(x_{ij})$, and $\I$ stands for the identity matrix.

A certain (di)graph may be a multidimensional circulant for several
different values of $n$. For instance, the digraph $G(\M;A)$ with
$\M=\mbox{ diag}(2,2,3)$  and  $A=\{ (1,0,0)^\top ,(0,1,0)^\top,$ $
(0,0,2)^\top \}$ is isomorphic to the digraph
$G(\M';A')$ with $\M'=\mbox{ diag}(2,6)$ and $A'=\{ (0,3)^\top ,(1,0)^\top ,
(0,4)^\top \}$ since, if $\U'$ is the $2\times 3$ matrix
obtained by taking the two last rows of the matrix $\U$ in (\ref{eq1}),
we have $\U'A\equiv A'\pmod{\M'}$. Following Leighton's terminology
\cite{r13}, if $k$ is the smallest value of such $n$ we will say that
the multidimensional circulant has {\it dimension} $k$ or that it is
{\it k-dimensional}. Notice that this parameter is in fact the minimum
rank of the groups such a (di)graph can arise from. Then the class
of circulants is precisely the class of $1$-dimensional circulants.

In studying the dimension of a given multidimensional circulant
$G(\M;A)$ we only need to consider the connected case. Indeed suppose
that the $\alpha$ disjoint components of $G(\M;A)$ are, say,
$k$-dimensional and isomorphic to $G(\M';A'), \M' \in \Z^{k \ast k}$. Then
$G(\M;A)$ is isomorphic to $G(\alpha \M';\alpha A')$ where $\alpha \M'$ and
$\alpha A'$ denote the matrix and set obtained from $\M'$ and $A'$ by
simply multiplying by ${\alpha }$ any, say the first, component of the
corresponding (column) vectors. As a corollary, the dimension of
$G(\M;A)$ cannot be greater than the cardinality of the minimum subset of
$A$  that  generates  $\Gamma=\langle \a_{1},\a_{2}, \ldots
,\a_{d}\rangle$.

To obtain other results about the dimension of multidimensional
circulants it is useful to introduce the concept of \'{A}d\'{a}m
isomorphism. Let $\M\in \Z^{n \ast n}$ and $\M'\in \Z^{n' \ast n'}$. Then
the multidimensional circulants $G(\M;A)$ and $G(\M';A')$ are said to be
{\it \'{A}d\'{a}m isomorphic} if there exists  an  isomorphism  $\phi$
between the groups $\Z^{n}/\M\Z^{n}$ and $\Z^{n'}/\M'\Z^{n'}$ such that
$\phi (A)=A'$. For instance, if $u$ is a unit of $\Z/m\Z$, that is
$\gcd (u,m)=1$, the circulants $G(m;A)$ and $G(m;uA)$ are \'{A}d\'{a}m
isomorphic. In \cite{r1} it was first conjectured that any two isomorphic
circulant digraphs are \'{A}d\'{a}m isomorphic, but in subsequent papers
more attention was paid to the corresponding statement for circulant
graphs. For instance, Djokovic \cite{r8} and Turner \cite{r17}
independently proved that \'{A}d\'{a}m's conjecture is true for circulant
graphs with prime order, and this is also the case for circulant
digraphs \cite{r11}. The first counter-examples to this conjecture,
both for graphs and digraphs were given by Elpas and Turner in
\cite{r11}. In \cite{r2}, Alspach and Parsons characterized, in terms of a
condition on automorphism groups, the validity of \'{A}d\'{a}m's
conjecture for a given order $m$. In particular, the authors used this
characterization to show that it holds for $m=p_{1}p_{2}$ where $p_{1}$
and $p_{2}$ are different primes. In \cite{r5}, Boesch and Tindell
conjectured that all isomorphic $2$-step circulant graphs are
\'{A}d\'{a}m isomorphic. This was independently proved in \cite{r7} and
\cite{r19}. The same result for $2$-step circulant digraphs was given in
\cite{r10}. In fact, Delorme, Favaron and Mah\'eo \cite{r7} proved some more
general results concerning Cayley (di)graphs of Abelian groups. Using
our terminology, they are stated in the following theorem.

\begin{theo}[Delorme et al. \cite{r7}]
\label{the3.1}
Let $\M\in \Z^{n \ast n}$ and $\M'\in \Z^{n' \ast n'}$ and suppose that
$A=\{ \a_{1},\a_{2} \}$ and
$A'=\{ \b_{1},\b_{2} \}$
are generating sets for $\Z^{n}/\M\Z^{n}$ and $\Z^{n'}/\M'\Z^{n'}$,
respectively. Then the two $($connected\/$)$ multidimensional circulant
digraphs $G(\M;A)$ and $G(\M';A')$ are isomorphic iff they are
\'{A}d\'{a}m isomorphic, except in the case when there exist two group
isomorphisms
\[ \phi : \Z^{n}/\M\Z^{n} \longrightarrow \Z/2\eta \Z\times \Z/2\Z,\quad
\eta \in \Z^{+} \]
and
\[ \phi ' : \Z^{n'}/\M'\Z^{n'} \longrightarrow \Z/4\eta \Z \]
such that $\phi(A)=\{ (1,0)^\top ,(1,1)^\top \}$ and
$\phi'(A')=\{1,2\eta +1\}$. $($In this case the two digraphs are isomorphic
but, clearly, they are not \'{A}d\'{a}m isomorphic.$)$ Moreover, this
result is also true for connected multidimensional circulant graphs if
we change $A$ and $A'$ by $\pm A$ and $\pm A'$, respectively.\ $\Box$
\end{theo}

Note that the exceptional case $\{ \Z/2\eta \Z\times \Z/2\Z$, $\a_{1}
=(1,0)^\top , \a_{2}=(1,1)^\top  \}$ could also be characterized
by  the  defining  relations  $\a_{1}+\a_{2}=\a_{2}+\a_{1}$ (Abelian group), $2\eta\a_{1}=\vec0$ and $2\a_{1}=
2\a_{2}$. Moreover, this last relation is equivalent to writing
${\rm o}(\a_{1}-\a_{2})=2$ (which holds indeed if we
substitute the above values in (\ref{eq9})).

Let $\H=\M\V$ be the Hermite normal form of the matrix $\M\in \Z^{n \ast n}$.
Let $\S=\mbox{ diag}(s_{1},s_{2}, \ldots ,s_{n})=\U\M\V$ be its Smith normal
form with $s_{1}=s_{2}= \cdots =s_{n-r}=1$ and consider the $r \times r$
and $r \times n$ matrices $\S'$ and $\U'$ defined as in Section 2. From
the results (\ref{eq3}), (\ref{eq4}), (\ref{eq6}) and (\ref{eq7}) given
there we have the following theorem.

\begin{theo}
\label{the3.2}
The multidimensional circulants $G(\M;A)$, $G(\H;A)$ and
$G(\S';\phi(A))$, where $\phi(A)=\{ \U'\a:\a\in A \}$, are \'{A}d\'{a}m
isomorphic.\ $\Box$
\end{theo}
As an example of application of this theorem, we can again consider the
two isomorphic multidimensional circulants with matrices $\M=\mbox{ diag}
(2,2,3)$ and $\M'=\mbox{ diag}(2,6)$ mentioned before.

\begin{coro}
\label{cor3.3}
Let $G(\M;A)$ be a $k$-dimensional circulant. Then $k \leq r$. In
particular, if $r=1$ $(d_{n-1}=s_{n-1}=1)$ such a $($di$)$graph is a
circulant.\ $\Box$
\end{coro}
From the above corollary and Proposition \ref{pro2.1}(b) we get

\begin{coro}
\label{cor3.4}
Let $G(\M;A)$ be a $k$-dimensional circulant with
$m=p_{1}^{r_{1}}p_{2}^{r_{2}} \cdots p_{t}^{r_{t}}$ vertices.
Then $k\leq \max\{r_{i}:1\leq i\leq t \}$. In particular, if $m$ is square free
(that is, $m$ is not divisible by the square of a prime) $G(\M;A)$ is a
circulant.\  $\Box$
\end{coro}
This coincides with the result obtained by Leighton in \cite{r13} for a
multidimensional circulant graph $G(\M;A)$ with $\M$ a diagonal matrix.

In the case of multidimensional $2$-step circulants we can give a
complete characterization of circulants and, hence, of their dimension.

\begin{theo}
\label{the3.5}
Let $\M$ be an $n\times n$ matrix with $(n-1)$th determinantal divisor
$d_{n-1}$. Let $A=\{ \a_{1},\a_{2} \}$ be a generating set of
$\Z^{n}/\M\Z^{n}$. Then the $($connected$)$ multidimensional $2$-step circulant
digraph $G(\M;A)$ $($respectively, graph $G(\M;\pm A)$$)$, on $m=|\det \M|$
vertices, is a circulant iff one of the following conditions holds:
\begin{itemize}
\item[$(a)$]
$d_{n-1} = 1$; or
\item[$(b)$] $d_{n-1} = 2$ and $m=2\gcd(m,\gcd(m\M^{-1}(\a_{1}-\a_{2})))$; $($respectively, or
\item[$(c)$]
$d_{n-1} = 2$ and
$m=2\gcd(m,\gcd(m\M^{-1}(\a_{1}+\a_{2})))$$)$.
\end{itemize}
\end{theo}
\noindent{\bf Proof.}
From Theorem \ref{the3.1} and Corollary \ref{cor3.3} it is
clear that $(a)$ is a necessary and sufficient condition for $G(\M;A)$ or
$G(\M;\pm A)$ to be a circulant except in the case $\Z^{n}/\M\Z^{n} \cong
\Z/2\eta \Z \times \Z/2\Z$ and o$(\a_{1}-\a_{2})=2$ (or possibly,
in the case of graphs, o$(\a_{1}+\a_{2})=2$.) (According to
this theorem, in this case we also have a circulant.) But then from the
results of Section 2, and in particular (\ref{eq8}), condition $(b)$
(or condition $(c)$, in the case of graphs) must hold.\ $\Box$

In \cite{r12,eaf} it was shown that the study of some distance-related
parameters, such as the diameter, of $2$-step circulant digraphs is best
accomplished by considering them as particular instances of
multidimensional circulants digraphs $G(\M;\e_{1},\e_{2})$
with $\M$ a $2\times 2$ matrix. Such (strongly connected) digraphs have
been called {\it commutative 2-step digraphs} \cite{r10}. The reason is
that the matrix $\M=(m_{ij})$ can always be chosen so that the studied
parameter is easily related to its entries $m_{ij}$. (In some cases the
same fact is true for $2$-step circulant graphs, see \cite{r4} and
\cite{r18}.) Hence, it is of some interest to characterize those
commutative $2$-step (di)graphs which are circulants. As a particular
case of Theorem 3.5, the next corollary gives a complete answer to this
question.

\begin{coro}
\label{cor3.6}
Let $\M=(m_{ij})$ be a $2\times 2$ integer matrix with $|\det \M|=m$. Then
the commutative $2$-step digraph $G(\M;\e_{1} ,\e_{2})$ is a circulant
digraph iff either
\[ d_{1}=\gcd (m_{11},m_{12},m_{21},m_{22})=1 \]
or
\[ d_{1}=2 \mbox{\ and\ } m=2\gcd (m,m_{22}+m_{12},m_{11}+m_{21}). \  \Box   \]
\end{coro}

Theorem \ref{the3.5}  illustrates the fact that, although the knowledge of the structure of
$\Z^{n}/\M\Z^{n}$ (Proposition \ref{pro2.1}) gives an upper bound for the dimension
of a multidimensional circulant (Corollary \ref{cor3.3}), the computation of its exact
value may require more sophisticated and particular techniques. This is also made
apparent for the next result, which gives the dimension of the direct product of
$n$ circulants, all of them with equal prime number of vertices. The  proof is similar to
that of Lemma 1 in \cite{r13}, which was suggested by Lawrence (personal
communication to Leighton.)

\begin{theo}
\label{the3.7}
Let $G_{1}=G(p;A_{1}), G_{2}=G(p;A_{2} ),\ldots ,G_{n}=G(p;A_{n})$ be $($connected\/$)$
circulants with $p>2$ vertices, $p$ a prime. Let $\M\in \Z^{n \ast n}$ be the diagonal
matrix {\rm diag}$(p,p,\ldots p)$ and $A=\{ a\e_{i}:a\in A_{i} ,1\leq i\leq n \}$. Then the
multidimensional circulant $G(\M;A)$ has dimension $n$.
\end{theo}

\noindent{\bf Proof.} First note that the (di)graph $G(\M;A)$, with vertex-set
$V=\Z/p\Z \times  \cdots  \times \Z/p\Z$ ($n$ factors), is nothing more than the
cartesian product $G_{1}\times \cdots \times G_{n}$. Let $\Omega$  be any
regular Abelian subgroup of ${\sl Aut}G(\M;A)$, $|\Omega |=|V|=p^{n}$. By Sabidussi's
result it suffices to show that $\Omega \cong \Z^{n}/\M\Z^{n} \cong
\Z/p\Z \times  \cdots  \times \Z/p\Z$.
For some fixed $1\leq i\leq n$ and $j\in \Z/p\Z$,  let $G_{ij}$  denote the sub(di)graph
of $G(\M;A)$ spanned by the vertices  whose labels have their $i$th component equal to $j$.
Now, let us consider the group $\Omega$  as acting on the set
${\cal G}=\{ G_{ij}:1\leq i\leq n, 0\leq j< p \}$.
To show that any automorphism $\omega \in \Omega$ preserves the set ${\cal G}$ (that is,
either $\omega(G_{ij}) \bigcap G_{kl} = \emptyset$ or $\omega(G_{ij}) = G_{kl}$),
it suffices to prove that $\omega$ preserves the set of $n$ directions. (As expected,  a
{\it direction} is defined as the set of edges whose endvertices only differ in a given
coordinate.) To know whether different edges belong to the same direction
we can apply the following algorithm:

\noindent Choose a vertex $\u \in V$ and consider any
shortest odd cycle containing it. Then, all the edges of this cycle clearly belong to one direction,
say $i$. If the cycle has length $p$, then the sub(di)graph spanned by its $p$ vertices,
$G_i(\u)$, is the copy of $G_{i}=G(p;A_{i})$  that contains vertex $\u$ and whose edges
belong to direction $i$.  Otherwise, we consider an edge of the cycle
and look for a
different shortest odd cycle (of the same length as before) containing it. In this way we
successively find the edges (and vertices) of $G_i(\u)$. If, in some step, there is no
such a shortest cycle we consider another of the (already found) edges of $G_i(\u)$.
Because of the nature of $G_{i}$, it is not difficult to realize that we eventually find the
searched $p$ vertices spanning $G_i(\u)$. To find the other sub(di)graphs $G_j(\u)$,
$j \neq i$, we start again from vertex $\u$ and look,  in the same way as before, for shortest
odd cycles not containing edges in the already found directions. This is done until no edge
incident to $\u$ is left out of discovered directions.
\noindent
To identify the directions of the
edges incident to other vertices, different from $\u$,  we can apply the
following procedure:
Consider two adjacent edges with endvertex $\u$ and different directions, say
$(\z,\u) \in G_j(\u)$ and $(\u,\v) \in G_i(\u)$. Look
for a shortest cycle containing them, $\u,\v,\ldots,\z,\u$,
(note that its length must be at least 4).
Then, the first edge (of the cycle) not in $G_i(\u)$ has direction $j$ and,
similarly, the last edge not in $G_j(\u)$  has direction $i$. Finally, once the
directions of a sufficient number of edges incident to a vertex, say $\w$, have been
determined, we can search again for appropriates odd cycles going through it,
in order  to locate the (di)graphs $G_i(\w)$, $i=1,2,\ldots ,n$.

From the above, the action of an automorphism  $\omega \in
\Omega$  on $V$ completely determines its action on ${\cal G}$. Conversely, let
$\u=(u_{1},u_{2}, \ldots ,u_{n}) \in V$. Then $\u$ is the only vertex
the sub(di)graphs $G_{iu_{i}}, 1\leq i\leq n$, have in common. Hence, the action of
$\omega$ on ${\cal G}$ also determines its action on $V$.

Let ${\cal G}_{1}, \ldots ,{\cal G}_{k}$ be the orbits of ${\cal G}$ under the action of
$\Omega$. Let  $\Omega_h$, $1\leq h\leq k$, be the restriction of  $\Omega$  to
${\cal G}_{h}$  with duplicates eliminated. Then,  $\Omega \subseteq  \Omega_1 \times
\cdots  \times \Omega_k$ and hence
\begin{equation}
\prod_{h=1}^k |\Omega_{h}| \geq |\Omega| = p^{n}.        \label{eq10}
\end{equation}
Moreover, since  $\Omega_h$  is Abelian and transitive on ${\cal G}_h$, it is also regular.
Therefore $|\Omega_h| = |{\cal G}_h|$, $1\leq h\leq k$, and then
\begin{equation}
\sum_{h=1}^k |\Omega_{h}| = \sum_{h=1}^k |{\cal G}_{h}| = |{\cal G}| = np.        \label{eq11}
\end{equation}
In addition, the order of an orbit, $|\Omega_h|= |{\cal G}_h|$, divides the order of the
permutation group $|\Omega |=p^n$, see \cite{r15}. Thus there exist integers $r_h \geq 0$,
$1\leq h\leq k$, such that $|\Omega_h|=p^{r_{h}}$ , and formulas (\ref {eq10}),
(\ref {eq11}) yield
\[
\sum_{h=1}^k r_h \geq n,\quad  \sum_{h=1}^k p^{r_{h}} = np,
\]
respectively. Hence, we must have $\sum_{h=1}^k pr_{h} \geq pn = \sum_{h=1}^k p^{r_{h}}$.
But, for $p>2$, $pr_{h} < p^{r_{h}}$ if $r_{h} \neq 1$ and $pr_{h}=p^{r_{h}}$ otherwise.
Thus $|\Omega_h|=p$ for any $1\leq h\leq k$, so that  $\Omega_h$  is isomorphic
to the cyclic group $\Z/p\Z$ and, from (\ref{eq10}),  $\Omega \cong \Z/p\Z \times  \cdots
\times \Z/p\Z$ ($n$ factors) as claimed. \ $\Box$

As an example of application of the above theorem,  we can state the following
special case:
\begin{coro}
\label{cor3.8}
The cartesian product  $K_p  \times  \cdots  \times K_p$ $($n factors$)$, and
the cartesian product of $n$ $p$-cycles $($directed or not$)$ both have dimension $n$.\ $\Box$
\end{coro}

In the above examples there is an easy way to know
whether edges incident to a given vertex $\u$ belong to the same direction (or to locate the
(di)graphs $G_i(\u)$), owing to the connected components of the neighbourhood of
$\u$ in the case of complete graphs $K_p$, and to the (shortest) $p$-cycles in the case of
cycles.
Another example comes when each set $A_i$ has the property that if $x, y \in A_i$,
$x \neq y$, then one at least of the elements $x-y$, $x+y$, $-x+y$, $-x-y$ belongs to $A_i$
(for example, if $A_i$ is stable under multiplication by 2): the directions are then given
by the connected components of the neigbourhood of $\u$.

\noindent
\large
{\bf Acknowledgment}

\normalsize
\noindent
This  work  has been supported in part by the Spanish Research Council (Comisi\' on
Interministerial de Ciencia y Tecnolog\'{\i}a, CICYT) under Projects TIC90-0712 and
TIC-92-1228-E.
The author thanks Charles Delorme for his valuable comments on the proof of
Theorem \ref{the3.7} and Corollary \ref{cor3.8},
and one of the referees for many helpful suggestions.



\begin{thebibliography}{99}

\bibitem{r1}
A. \'{A}d\'{a}m, Research problem 2-10, {\it J. Combin. Theory}
{\bf 2} (1967) 393.

\bibitem{r2}
B. Alspach and T.D. Parsons, Isomorphism of circulant graphs and
digraphs, {\it Discrete Math.} {\bf 25} (1979) 97--108.

\bibitem{r3}
J.-C. Bermond, F. Comellas and D.F. Hsu, Distributed loop computer
networks: a survey, to appear in {\it J. Parallel and Distributed
Computing} {\bf 24} (1995) 2--10.

\bibitem{r4}
J.-C. Bermond, G. Illiades and C. Peyrat, An optimization problem in
distributed loop computer networks, Proc. 3rd Int. Conf. Combin. Math.,
New  York, 1985. {\it Ann. of the New  York Acad. Sci.}
{\bf 555} (1989) 45--55.

\bibitem{r5}
F. Boesch and R. Tindell, Circulants and their connectivities,
{\it J. Graph Theory} {\bf 8} (1984) 487--499.

\bibitem{r6}
G. Chartrand and L. Lesniak, {\it Graphs and Digraphs},
Wadsworth, Monterrey, 1986.

\bibitem{r7}
C. Delorme, O. Favaron and M. Mah\'eo, Isomorphisms of Cayley
multigraphs of degree 4 on finite Abelian groups,
{\it European J. Combin.} {\bf 13} (1992), {\it no.} 1, 59--61.

\bibitem{r8}
D.Z. Djokovic, Isomorphism problem for a special class of
graphs, {\it Acta Math. Acad. Sci. Hung.} {\bf 21} (1970) 267--270.

\bibitem{r11}
B. Elpas and J. Turner, Graphs with circulant adjacency matrices,
{\it J. Combin. Theory} {\bf 9} (1970) 297--307.

\bibitem{eaf}
P. Esqu\'e, F. Aguil\'o and M.A. Fiol, Double commutative-step digraphs
with minimum diameters,
{\it Discrete Math.} {\bf 114} (1993) 147--157.

\bibitem{r9}
M.A. Fiol, Congruences in $Z^n$, finite Abelian groups and the Chinese
remainder theorem, {\it Discrete Math.} {\bf 67} (1987) 101--105.

\bibitem{r10}
M.A. Fiol and P. Morillo, Congruences in $Z^2$ and commutative
two step digraphs, {\it 3\`{e}me Coll. Int. Theorie des Graphes et
Combinatoire}, Marseille, France, June 1986.

\bibitem{r12}
M.A. Fiol, J.L.A. Yebra, I. Alegre and M. Valero, A discrete
optimization problem in local networks and data alignment,
{\it IEEE Trans. Comput.} {\bf C-36} (1987) 702--713.

\bibitem{r13}
F.T. Leighton, Circulants and the characterization of
vertex-transitive graphs, {\it J. Res. Natl. Bur. Standards} {\bf 88}
(1983), no. 6, 395--402.

\bibitem{r14}
M. Newman, {\it Integral Matrices}, Pure and Appl. Math. Series
Vol. {\bf 45}, Academic Press, New  York, 1972.

\bibitem{r15}
D.J.S. Robinson, {\it A Course in the Theory of Groups}, Springer, New  York, 1982.

\bibitem{r16}
G. Sabidussi, On a class of fixed-point-free graphs, {\it Proc.
Amer. Math. Soc.} {\bf 9} (1958) 800--804.

\bibitem{r17}
J. Turner, Point-symmetric graphs  with a prime number of points,
{\it J. Combin. Theory} {\bf 3} (1967) 136--145.

\bibitem{r18}
J.L.A. Yebra, M.A. Fiol, P. Morillo and I. Alegre,
The diameter of undirected graphs associated to plane tessellations,
{\it Ars Combin.} {\bf 20B} (1985) 159--171.

\bibitem{r19}
S.C. Zhou, On \'{A}d\'{a}m isomorphism of circulants, {\it J. Changsha
Railway Inst.} {\bf 5} (1987), no. 2, 11--18.

\end{thebibliography}
\end{document}